\documentclass{elsart3-1}

\usepackage{setspace}

\usepackage{amsmath,amssymb}
\usepackage{color,xcolor}
\usepackage{wrapfig}  
\usepackage{graphicx}
\usepackage{enumerate}
\usepackage{cases}
\usepackage{verbatim}
\renewcommand{\H}{\textbf{H}}

\newcommand{\x}{\textbf{x}}

\newcommand{\E}{\textbf{E}}
\newcommand{\J}{\textbf{J}}

\newcommand{\gc}{\textbf{g}_C}
\newcommand{\gcnm}{\textbf{g}_{C}^{n,m}}
\newcommand{\gi}{\textbf{g}_I}
\newcommand{\ginm}{\textbf{g}_{I}^{n,m}}

\renewcommand{\Re}{\mathrm{Re}\,}
\renewcommand{\Im}{\mathrm{Im}\,}

\newcommand{\Om}{\Omega}
\newcommand{\Gm}{\Gamma}
\newcommand{\bnu}{{\bf n}}

\newcommand{\Tau}{\Gm}

\newcommand{\bmu}{\mu}

\newcommand{\bp}{p}

\newcommand{\Mnm}{{\bf M}_{n}^{m}}

\newcommand{\Ynm}{{\bf Y}_{n}^{m}}
\newcommand{\Unm}{\textbf{U}_{n}^{m}}
\newcommand{\Vnm}{\textbf{V}_{n}^{m}}
\newcommand{\anm}{{a}_{n}^{m}}

\newcommand{\bnm}{{b}_{n}^{m}}

\newcommand{\cnm}{{c}_{n}^{m}}
\newcommand{\dnm}{{d}_{n}^{m}}

\DeclareMathOperator{\curl}{\textbf{curl}}
\DeclareMathOperator{\curlS}{\text{curl}_{\Tau}}
\DeclareMathOperator{\CurlS}{\overrightarrow{\text{curl}}_{\Tau}}

\newcommand{\grad}{\nabla}
\DeclareMathOperator{\gradS}{\nabla_{\Tau}}
\DeclareMathOperator{\gradSo}{\nabla_{\ci}}

\DeclareMathOperator{\dvgS}{\text{div}_{\Tau}\hspace{.01in}}

\newcommand{\sumnm}{\sum_{n=1}^{\infty}\sum_{m=-n}^{n}}

\newcommand{\ci}{{\mathbb{S}^2}} 

\usepackage{amssymb}

\usepackage[english,francais]{babel}
\usepackage[latin1]{inputenc}
\newtheorem{theorem}{Theorem}[section]
\newtheorem{lemma}[theorem]{Lemma}
\newtheorem{e-proposition}[theorem]{Proposition}

\newtheorem{e-definition}[theorem]{Definition\rm}

\newtheorem{theoreme}{Th\'eor\`eme}[section]

\newtheorem{proposition}[theoreme]{Proposition}

\renewcommand{\theequation}{\arabic{equation}}
\def\og{\leavevmode\raise.3ex\hbox{$\scriptscriptstyle\langle\!\langle$~}}
\def\fg{\leavevmode\raise.3ex\hbox{~$\!\scriptscriptstyle\,\rangle\!\rangle$}}

\journal{the Acad\'emie des sciences}
\begin{document}
\centerline{}
\begin{frontmatter}


\selectlanguage{english}

\title{Improved non-overlapping domain decomposition algorithms for the eddy current problem}
  

\selectlanguage{english}



%

\author[NJIT]{Y. Boubendir},
\ead{boubendi@njit.edu}
\author[CMAP]{H. Haddar}
\ead{haddar@cmap.polytechnique.fr}
\author[NJIT]{M. K. Riahi}
\ead{riahi@njit.edu}
\ead[url]{https://web.njit.edu/~riahi/}
\address[NJIT]{Department of mathematical science, New Jersey Institute of Technology, New Jersey USA.}
\address[CMAP]{INRIA,  Ecole Polytechnique (CMAP) and Universit\'e Paris Saclay, Route
  de Saclay, Palaiseau 91128 Cedex FRANCE.}



\medskip
\begin{center}
{\small Received *****; accepted after revision +++++\\ Presented by Â£Â£Â£Â£Â£}
\end{center}

\begin{abstract}
A domain decomposition method is proposed based on carefully chosen impedance
transmission operators for a hybrid formulation of the eddy current
problem. Preliminary analysis and numerical results are provided in the
spherical case  showing the potential of these conditions in accelerating the
convergence rate.
\selectlanguage{english}
{\it To cite this article:
Y. Boubendir, H. Haddar and M.K. Riahi, C. R. Acad. Sci. Paris, Ser. I 340 (2015).}

\vskip 0.5\baselineskip

\selectlanguage{francais}
\noindent{\bf R\'esum\'e} \vskip 0.5\baselineskip \noindent
Nous proposons une méthode de décomposition de domaine basée sur un choix
particulier de l'écriture des conditions de transmission pour une formulation hybride
du modèle courant de Foucault. Nous donnons des résultats analytiques et
numériques préliminaires dans le cas sphérique qui montrent le potentiel de
ces conditions dans l'accélération de la vitesse de convergence d'une résolution
itérative du problème.  
{\it Pour citer cet article~: Y. Boubendir, H. Haddar et M.K. Riahi, C. R. Acad. Sci.
Paris, Ser. I 340 (2015).}

\end{abstract}
\end{frontmatter}

\selectlanguage{francais}
\section*{Version fran\c{c}aise abr\'eg\'ee}
Le modèle de courant de Foucault s'écrit pour le champ électrique $\E$ et le champ magnétique $\H$ 
\begin{equation}\label{eddyGlobal}
\begin{cases}
\curl \H - \sigma \E  &=0  \quad\text{ dans }\Omega \\
\curl \E - i\omega\bmu\H&=0\quad\text{ dans }\Omega\\
\bnu\times \H &=\bnu\times\J_e\quad\text{ sur }\partial\Omega.
\end{cases}
\end{equation}
où la perméabilité magnétique $\mu$ est une fonction à valeurs réelles qui peut
dépendre de l'espace, et $\J_e$ représente un terme source. La conductivité
électrique $\sigma$ est nulle en dehors de la partie conductrice $\Om_C$. Ces equations sont
complétées par les conditions d'interface traduisant la continuité des
composantes tangentielles de $\E$ et $\H$ à travers $\Gamma$, l'interface entre
$\Om_C$ et le vide $\Om_I:=\Om\backslash \overline \Om_C$. 

Le champ magnétique $\H$ est à rotationnel nul dans $\Om_I$.  Ceci implique,
lorsque ce dernier est simplement connexe, que $\H=\grad \bp$ (voir par exemple
\cite{MR2680968}), et le modèle s'écrit sous la forme
\begin{equation}\label{Orgproblem}
\begin{cases}
\curl\curl \E - \kappa^2\E  = 0 \quad \text{dans}\quad\Omega_C \\
\nabla \cdot (\mu \nabla p) = 0 \quad \text{dans}\quad\Omega_I
\end{cases}
\end{equation}
avec $\kappa^2:=i\omega\mu\sigma$. Les conditions de continuité à l'interface $\Gamma$ s'écrivent
\begin{equation}\label{orgcondts}
\curl\E \times \bnu =  i\omega\mu\CurlS p \quad \quad 
\frac{\partial p}{\partial\bnu} = \frac{1}{i\omega\mu}\curlS\E
\quad \text{sur}\quad\Gamma
\end{equation}
où $\bnu $ est un vecteur normal unitaire sur $\Gamma$ qui pointe vers
l'extérieur de $\Om_C$ et où $\CurlS p:= \gradS p
\times \bnu$ et $\curlS\E := \dvgS(\E \times \bnu)$ avec $\gradS$ et $\dvgS$
désignant respectivement le gradient et la divergence surfaciques. Nous
remarquons que le cas non simplement connexe pourrait également être traité au
moyen de rajouter à $p$ les contributions d'un nombre fini de fonctions
harmoniques (\cite{MR2680968}).
Dans ce travail nous proposons une méthode de décomposition de domaine basé sur
la décomposition de $\Omega$ en partie conductrice et le vide mais en
remplaçant les conditions d'interface naturelle par les conditions de continuité
suivantes 
\begin{eqnarray}
    \frac{\partial p}{\partial\bnu} + \beta_I \curlS\CurlS p = \frac{1}{i\omega\mu}\curlS \left(\E+\beta_I \curl\E\times\bnu\right) \label{BCI}\\
   \curl\E\times\bnu + \beta_C\CurlS \left( \curl\E\cdot\bnu\right)= i\omega\mu\CurlS\left( p + \beta_C    \frac{\partial p}{\partial\bnu} \right) \label{BCC}
\end{eqnarray}
où $\beta_I $ and $\beta_C$ sont des nombres complex appropriées. On montre que
ces conditions sont équivalentes aux conditions originales sous certaines
conditions sur $\beta_C$ et $\beta_I $ (Lemme \ref{consistance}). Notre
algorithme itératif est décrit par \eqref{lap}-\eqref{Max}. Dans la cas sphérique, on montre
que asymptotiquement, l'opérateur d'itération \eqref{defT} est une contraction plus rapide
que celle correspondant au choix $\beta_I= \beta_C =0$ (Proposition
\ref{asymptotics}). Nous concluons cette note par un exemple numérique $3D$ montrant une
meilleure vitesse de convergence pour des valeurs de $\beta_C $ et $ \beta_I$
différentes de $0$.

\selectlanguage{english}
\section{Introduction}\label{Intro}
The eddy current approximation of the Maxwell equations, for the electric field $\E$ and the magnetic field $\H$ reads
\begin{equation}\label{eddyGlobal}
\begin{cases}
\curl \H - \sigma \E  &=0  \quad\text{ in }\Omega \\
\curl \E - i\omega\bmu\H&=0\quad\text{ in }\Omega\\
\bnu\times \H &=\bnu\times\J_e\quad\text{ on }\partial\Omega.
\end{cases}
\end{equation}
where the magnetic permeability $\mu$ is a real valued function that
may depend on space and  $\J_e$ stands for the source excitations.
The electric conductivity $\sigma$ has support only in the conductive 
material $\Om_C$. These equations are supplemented with the continuity of
tangential components of $\E$ and $\H$ across $\Gamma$, the interface between
conductive and non conductive regions.

The field $\H$ is 
curl free in the insulating region $\Om_I:=\Om\backslash \overline \Om_C$. 
 This implies in the
case of a simply connected topology that  the magnetic field is a
gradient of a harmonic scalar field  i.e. $\H=\grad \bp$ (see~
\cite{MR2680968} and references therein), which leads 
 to the problem
\begin{equation}\label{Orgproblem}
\begin{cases}
\curl\curl \E - \kappa^2\E  = 0 \quad \text{in}\quad\Omega_C \\
\nabla \cdot (\mu \nabla p) = 0 \quad \text{in}\quad\Omega_I.
\end{cases}
\end{equation}
where $\kappa^2:=i\omega\mu\sigma$. The interface continuity conditions across $\Gamma$ lead to 
\begin{equation}\label{orgcondts}
\curl\E \times \bnu =  i\omega\mu\CurlS p \quad \quad 
\frac{\partial p}{\partial\bnu} = \frac{1}{i\omega\mu}\curlS\E
\quad \text{on}\quad\Gamma
\end{equation}
where $\bnu $ is a unitary normal vector on $\Gamma$ pointing toward the exterior of the conductor. We use $\CurlS p:= \gradS p
\times \bnu$ and $\curlS\E := \dvgS(\E \times \bnu)$ with $\gradS$ and $\dvgS$
respectively being the surface gradient and the surface divergence. To simplify
the presentation we also assume that $\partial \Omega \subset \partial
\Omega_I$ and close the problem by imposing a boundary condition 
$ p = f \text{ on }\partial \Omega.$ Let us already remark that the case of non
simply connected domain can also be treated in a similar way but with 
additional technicalities related to the incorporation of (finite number of) divergence and curl
free functions.

Many methods, such as potential and hybrid formulations
\cite{MR3349682,MR3343738,NUMNUM20060},  have been
developed
to deal with this type of problem. An extensive overview of these methods can
be found in \cite{MR2680968}. 
Generally speaking, these methods fall in the class of  direct
formulations gathering the problem in $\Omega_C$ and $\Omega_I$ into the same
linear system to solve. A nice discussion of the advantages/disadvantages of such methods and others
can be also found in \cite{MR3343738}, where the authors propose 
a new numerical discretization to overcome these difficulties, mainly related 
to the size of the linear systems. Although the  approach  in \cite{MR3343738}
represents an important contribution in solving these difficulties, 
it remains strongly dependent on the use of preconditioners because of
the ill-conditioning which is further worsen by the high contrast created by the 
conductivity in $\Omega_C$.

Domain decomposition methods~\cite{quarteroni1999domain,toselli2005domain} are well suited 
for this problem since they allow the problems 
in $\Omega_C$ and $\Omega_I$ to be solved separately with appropriate approaches.
Up to our knowledge, the only domain decomposition algorithm developed for formulation \eqref{Orgproblem} is
the given in \cite{MR2680968} and exploits transmission conditions \eqref{orgcondts}. 
We  here propose to improve this method by modifying these conditions under the form 
\begin{eqnarray}
    \frac{\partial p}{\partial\bnu} + \beta_I \curlS\CurlS p = \frac{1}{i\omega\mu}\curlS \left(\E+\beta_I \curl\E\times\bnu\right) \label{BCI}\\
   \curl\E\times\bnu + \beta_C\CurlS \left( \curl\E\cdot\bnu\right)= i\omega\mu\CurlS\left( p + \beta_C    \frac{\partial p}{\partial\bnu} \right) \label{BCC}
\end{eqnarray}
where $\beta_I $ and $\beta_C$ denote appropriate (complex) numbers. For
similar ideas in different contexts we refer the reader to~\cite{boubendir2007analysis,gander2002optimized,boubendir2012quasi}.

We establish in this part a consistency property for the above impedance boundary condition.

\begin{lemma}\label{consistance}
The following conditions 
\begin{equation}\label{condProd}
\Re\{-\beta_C\beta_I\} \geq 0, \quad\text{ or }\quad \Im\{\beta_C\beta_I\}\neq 0
\end{equation} ensure consistency
between  the original conditions \eqref{orgcondts} and the new ones
\eqref{BCI}-\eqref{BCC}. 
\end{lemma}
\textit{Proof.}
Let us define the following quantities on $\Gamma$
\begin{equation}
\Xi_{\text{C}}:=\curl\E \times \bnu - i\omega\mu\CurlS\bp, \qquad \Xi_{\text{I}}:=\dfrac{\partial\bp}{\partial\bnu} - \dfrac{1}{i\omega\mu}\curlS\E,
\end{equation}
which are zero if interface conditions \eqref{orgcondts} 
 hold. Interface conditions \eqref{BCC} and \eqref{BCI} can be written as 
\begin{equation}
 \Xi_{\text{C}} + i\omega\mu \beta_C\CurlS \Xi_{\text{I}} =0 ,\qquad
 \Xi_{\text{I}} + \dfrac{1}{i\omega\mu} \beta_I\curlS \Xi_{\text{C}} =0.
\end{equation}
 Therefore, $\Xi_{\text{C}} - \beta_C \beta_I \CurlS\curlS \Xi_{\text{C}} =0.$
with implies, using a variational form, that $\|\Xi_{\text{C}}\|^2_{L^2(\Gamma)} - \beta_C \beta_I \| \curlS \Xi_{\text{C}}\|^2_{L^2(\Gamma)}  =0$.
If \eqref{condProd} is satisfied then $\Xi_{\text{C}}=0$ which also implies
$\Xi_{\text{I}} =0$ and one recovers  interface conditions \eqref{orgcondts}.
$\hfill\square$

	 The problem at hand is then reduced to the following iterative algorithm where the two problems are solved separately (with appropriate variational formulations and discretization) 
\begin{eqnarray}
&&\left\{\begin{array}{lll}
\Delta p^{(k+1)} = 0 \quad \text{in}\ \Omega_I \\
    \dfrac{\partial p^{(k+1)}}{\partial\bnu} + \beta_I \curlS\CurlS p^{(k+1)} &=\dfrac{1}{i\omega\mu}\curlS  (\E^{(k)}+\beta_I \curl\E^{(k)}\times\bnu ) \,\,&\text{on}\,\Gamma \\
\end{array}\right.\label{lap} \\
&&\left\{\begin{array}{lll}
\curl\curl \E^{(k+1)} - \kappa^2 \E^{(k+1)}= 0 \quad \text{in}\ \Omega_C\\
\curl\E^{(k+1)}\times\bnu + \beta_C\CurlS ( \curl\E^{(k+1)}\cdot\bnu )&=  i\omega\mu\CurlS ( p^{(k)} + \beta_C    \dfrac{\partial p^{(k)}}{\partial\bnu}  ) &\text{on}\,\Gamma.
\end{array}\right.\label{Max}
\end{eqnarray}
It order to ensure well posed problems for $p^{(k+1)}$ and $\E^{(k+1)}$, a
  variational study of problems \eqref{lap} and \eqref{Max} show that sufficient conditions are 
respectively $\Re \beta_I \leq 0$ and $\Re \beta_C \geq 0, \; \Im \beta_C \leq 0$.\vspace{-.5cm}

\section{Iteration operator in the case of concentric spheres}\vspace{-.2cm}
This section is dedicated to the computation of the eigenvalues of
the iteration operator, denoted by $\mathcal{T}$,  
in the case of concentric spheres. 
The main goal is to study the dependence of these eigenvalues on $\beta_I$ and $\beta_C$ 
 and show that non zero values of $\beta_I$ and $\beta_C$ improve the
 convergence of the iterative procedure. In the case $f=0$ (meaning the
 solution is $0$) better behavior correspond with eigenvalues of modulus closer to $0$ in the case of a successive iterative procedure.
Consider the sphere $\Om=B_{R}\subset\mathbb R^3$ with
radius $R>1$. The insulating and conducting regions are  respectively
given by $B_{R} \backslash B_{1}$ and $B_{1}$, where $B_{1}$ 
is the unit sphere (the case of a sphere of radius $r <R$ can be easily deduced
using an appropriate scaling).  Assume that $\gc$ (resp. $\gi$)
represents the right side on $\Gamma$ in \eqref{Max} (resp. \eqref{lap}), and let us define  
$\mathbf{g}=(\gc,\gi)^T$. Performing one iteration consists in computing
\begin{equation} \footnotesize\label{defT}
 \mathbf{g}^{(n+1)}= \mathcal{T} \mathbf{g}^{(n)} =\left(\begin{array}{cc}
      0 & \;\;  \mathcal{T}_C \\[4pt]
      \mathcal{T}_I & \;\; 0
    \end{array}\right)\left(\begin{array}{cc}
      \gc^{(n)} \\[8pt]
      \gi^{(n)}
    \end{array}\right)
\end{equation}
with
\begin{equation}\label{Equ10}
 \mathcal{T}_I \left(\gc^{(n)} \right):=\dfrac{1}{i\omega\mu}\curlS \left(\E^{(n+1)}
 +\beta_I \curl\E^{(n+1)}\times\bnu\right), \quad 
\mathcal{T}_C  \left(\gi^{(n)} \right) := i\omega\mu\CurlS\left( p^{(n+1)} + 
\beta_C \partial_{\bnu} p^{(n+1)} \right)
\end{equation}
 Let  $\Ynm$, $n=0,1, \ldots$, $ -n \le m \le n$ denote the spherical harmonics
 and set
 $\Unm := \dfrac{1}{\sqrt{n(n+1)}}\gradSo\Ynm$ 
  $\text{and}\quad \Vnm(\hat \x) := \hat \x\times\Unm(\hat \x)$, $\hat \x \in \ci$. We remark from the expression of
   $\gc$ that this field belongs to span $\{\Vnm, n=1,2 \ldots,  -n \le m \le n
   \}$. We also denote by $j_n$ the spherical Bessel function of the first kind
   of order $n$.
\begin{proposition}\label{eigenvalues}
 Assume that 
$$\gc=\sumnm \gcnm \Vnm \mbox{ and } \gi=\sumnm\ginm \Ynm(\hat\x).$$
 Then we have 
$$\mathcal{T}_C(\gi)=\sumnm (\mathcal{T}_C)^{m}_{n} \ginm \Vnm \mbox{ and } \mathcal{T}_I(\gc)=\sumnm (\mathcal{T}_I)^{m}_{n} \gcnm\Ynm(\hat\x)$$
with
 \begin{eqnarray}
(\mathcal{T}_C)^{m}_{n} :=\textstyle -i\omega\mu\left( \dfrac{1}{\sqrt{n(n+1)}}\mathbb{B}_I + \beta_C{{n(n+1)}}\mathbb{A}_I\right) \slash \left(\mathbb A_I + \beta_I \mathbb{B}_I \right) \label{tcgcnm} \\
(\mathcal{T}_I)^{m}_{n}  :=  \textstyle \dfrac{1}{i\omega\mu}\left(
  \dfrac{-1}{\sqrt{n(n+1)}} \mathbb{B}_C+\beta_I {\sqrt{n(n+1)}}\mathbb{A}_C \right) \Big\slash \left( \mathbb{A}_C+\beta_C\mathbb{B}_C \right)\label{tigcnm}
\end{eqnarray}
and $\mathbb A_I := -n\left( 1 + R^{-2n} \right ) $ , $\mathbb{B}_I := n(n+1)\left(1 - R^{-2n}\right)$,  $\mathbb{A}_C:=-\sqrt{n(n+1)}  ( j^{}_{n}(\kappa) + \kappa j^{'}_{n} (\kappa)  )$ and $\mathbb{B}_C:=- (n(n+1) )^{\frac 3 2} j^{}_{n}(\kappa)$.
 \end{proposition}
\textit{Proof.} 
The solution $p$ of problem \eqref{lap} can be expressed as
$
\bp(\x) = \sumnm \left( \cnm |\x|^n + \dnm |\x|^{-n}  \right) \Ynm (\hat \x).
$
From the homogeneous exterior boundary conditions we get  $\cnm = -\dnm  R^{-2n}$.
Using the boundary condition \eqref{BCI} we obtain 
 $$
  \left(\mathbb A_I + \beta_I \mathbb{B}_I \right)\dnm =  
  \ginm.
 $$ 
 The expression of $(\mathcal{T}_C)^{m}_{n}$ then follows from
 evaluating the expression  of $\mathcal{T}_C(\gi)$.

The solution $\E$ of problem \eqref{Max} can be expressed as 
$
\E (\x) =  \sumnm \anm \Mnm (\x) - \dfrac{i}{\kappa}\bnm \curl \Mnm(\x)
$
with $\Mnm(\x):=\curl\big(\x j_{n}(\kappa |\x|)\Ynm(\hat \x)\big)$. Since 
\begin{eqnarray*}
\Mnm(\x) =- \sqrt{n(n+1)}j^{}_n(\kappa |\x|)  \Vnm(\hat \x) \mbox{ and }
\curl\Mnm(\x) =- \sqrt{n(n+1)}\left( \dfrac {1}{|\x|} j^{}_n(\kappa |\x|) + \kappa j^{'}_n(\kappa |\x|) \right) \Unm(\hat \x) 
\end{eqnarray*}
we get at $|x|=1$,
\begin{eqnarray}\label{curlExn}
\curl\E\times\bnu =& \sumnm &-\anm \sqrt{n(n+1)} \left(
  j^{}_{n}(\kappa) + \kappa j^{'}_{n} (\kappa) \right) \Vnm(\hat \x) \notag+ i\kappa \bnm \sqrt{n(n+1)} j^{}_{n}(\kappa) \Unm(\hat \x).
\end{eqnarray}
In addition, because $\curlS\Unm=0$, we obtain at  $|x|=1$,
\begin{eqnarray}\label{CurlSurlS2E}
\CurlS\curlS \E = \sumnm -\anm \left(n(n+1)\right)^\frac 3 2 j^{}_{n}(\kappa) \Vnm(\hat \x).
\end{eqnarray}
Therefore, combination of \eqref{curlExn},\eqref{CurlSurlS2E} and  the boundary condition \eqref{BCC} leads to  
\begin{eqnarray}
\curl\E\times\bnu + \beta_C \CurlS\curlS \E = \sumnm \left( \mathbb{A}_C + \beta_C \mathbb{B}_C \right) \Vnm(\hat \x).
\end{eqnarray}
This implies in particular that   
\begin{equation}
 \anm= \gcnm \slash \left( \mathbb{A}_C+\beta_C\mathbb{B}_C \right), \quad\bnm=0.
\end{equation}
The expression $(\mathcal{T})^{m}_{n}$ then directly follows from evaluating $\mathcal{T}_I(\gc)$.
 $\hfill\square$

Using the structure of $\mathcal{T}$  one observes that the operator $\mathcal{T}^2$ is diagonal with eigenvalues\begin{equation}\label{amplcoef}
(\mathcal{T})^{m}_{n}:= (\mathcal{T}_I)^{m}_{n} \cdot (\mathcal{T}_C)^{m}_{n}.
\end{equation}
\begin{proposition}\label{asymptotics}
The leading term in the asymptotic expansion of the amplification coefficient for large $n$ and any
$\beta_I, \beta_C$ is given by
$$|(\mathcal{T})^{m}_{n}| \sim
\left| \frac{1-n \beta_C}{1+n\beta_C}\right| \quad \text{ as } n\rightarrow\infty.$$
\end{proposition}
\textit{Proof.}
	The asymptotic expansion for the Bessel functions of the first kind with
        complex argument~\cite[Formula 9.3.1]{abramowitz1964handbook} for a
        fixed complex argument $\kappa$ and a large integer $n$ is given by 
$ j_{n}(\kappa) \sim \frac{1}{\sqrt{2\pi n}}
\left(\frac{e\kappa}{2n}\right)^{n}$ and $ j_{n}'(\kappa) \sim
\frac{n}{\kappa}j_{n}(\kappa)$.We then have the following asymptotic expansions of the coefficients appearing in Proposition~\ref{eigenvalues} 
$\mathbb{A}_I\sim - {n}, \mathbb{B}_I\sim {n^2}$,
$\mathbb{A}_C\sim-{n^2} j_n(\kappa)$ and
$\mathbb{B}_C\sim-{n^3}j_n(\kappa)$. Plugging these asymptotics in formula
\eqref{tcgcnm} and \eqref{tigcnm}  respectively give  for large $n$
$$ (\mathcal{T}_C)^{m}_{n}\sim i\omega\mu ( 1- n\beta_C )/( 1-n\beta_I ) \mbox{
  and } (\mathcal{T}_I)^{m}_{n}\sim \frac{1}{i\omega\mu} (1 - \beta_I n )\slash( -1 - \beta_C  n).$$
This gives the announced result.  $\hfill\square$

We then conclude that as long as $\Re \beta_C >0$ the coefficient
$(\mathcal{T})^{m}_{n}$ has a modulus strictly smaller than $1$ for large
$n$. It is surprising that the asymptotic behavior is independent from the
values of $\beta_I$ which may suggest that this parameter has less important
influence in accelerating the iterations. More importantly, the asymptotic
formula shows that $|(\mathcal{T})^{m}_{n}( \Re \beta_C >0)| <
|(\mathcal{T})^{m}_{n}( \Re \beta_C =0)|$ for $n$ sufficiently large which
would indicate  for the cases  $\Re \beta_C >0$ better convergence properties than for
 the natural choice $\beta_I=0$ and $\beta_C=0$. We also remark that purely
 imaginary values of $\beta_C$ do not improve the convergence properties
 for large modes.

\begin{figure}[!hbp]
\centering
\includegraphics[height=5cm,width=5cm]{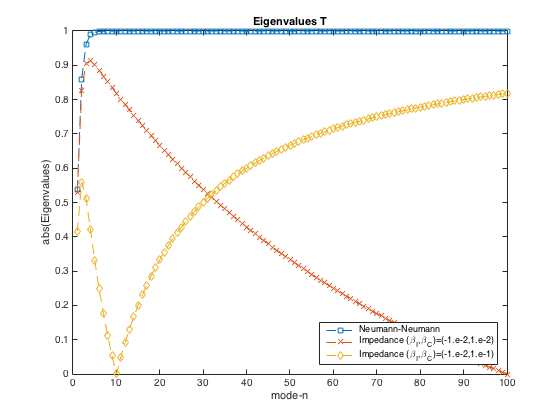}
\includegraphics[height=5cm,width=5cm]{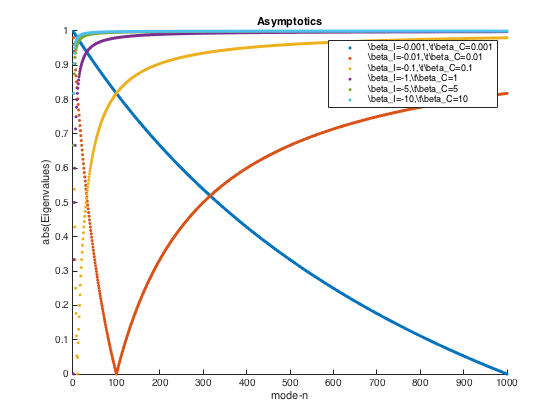}
\caption{Plots of $n\mapsto |(\mathcal{T})^{m}_{n}|$ for the choice of ($\beta_I,\beta_C$)=\{(-1.e-2,1.e-2),(-1.e-2,1.e-1) (left) and  different choices of $\beta_I=-\beta_C$  (right)\}.}\label{eigenplot}
\end{figure}

We plot in Figure~\ref{eigenplot} the amplification
coefficient~\eqref{amplcoef} of the iterative procedure for each mode $n,
m$. We observe that the above conclusions also hold for all modes and that
better behavior is observed for $\Re \beta_C >0$. The asymptotic behavior of
$|(\mathcal{T})^{m}_{n}|$ is also confirmed by Figure~\ref{eigenplot} right. 

\section{3D Finite Elements preliminary experiments}

\begin{figure}[htbp]\centering
\begin{tabular}{cccc}
&\includegraphics[height=5cm,width=5cm]{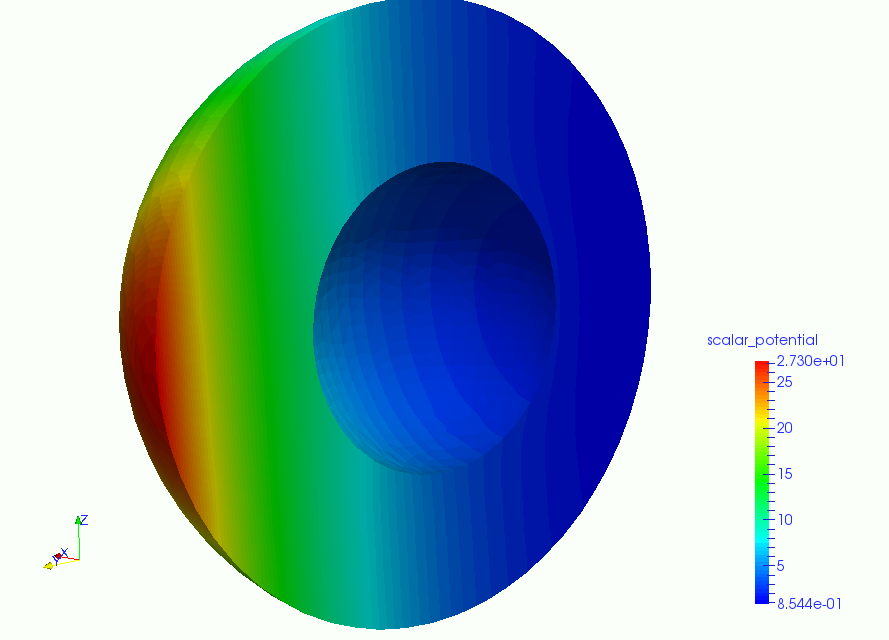}
&\includegraphics[height=5cm,width=5cm]{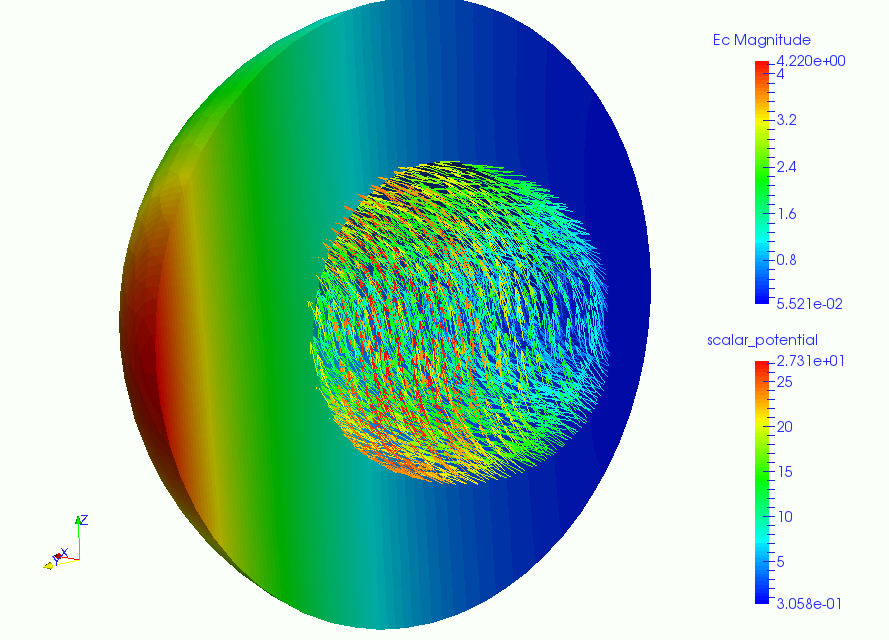}
&\includegraphics[height=5cm,width=5cm]{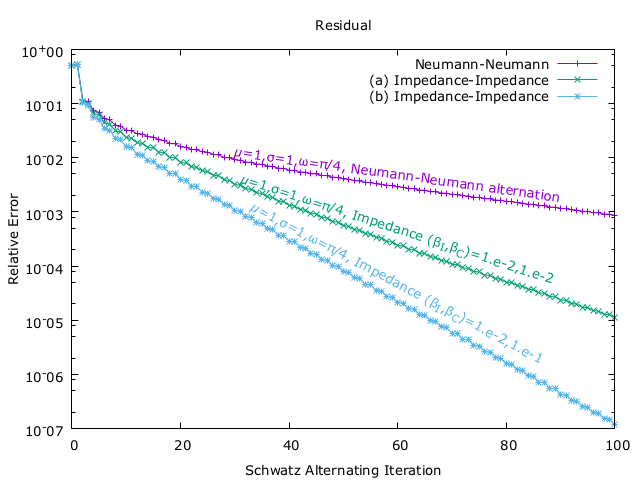}\\
&(a) &(b) &(c)
\end{tabular}
\caption{Computed solution for the potential $\bp$ in the insulator (a) and for
  electric field  in the conductor (b). The residual of the DDM iterations (c).}\label{presults}
\end{figure}

Preliminary results of the new formulation are presented in 
this section. The insulating region is given by $B_{R} \backslash B_{r}$ where $R=2$ and $r=1$, and the conduction region by
$B_r$. The electromagnetic coefficients are chosen to be; $\sigma=1$
for the conductivity and $\mu=1$ for the permeability. The frequency 
is set to $\omega=\pi/4$. The coefficients $\beta_I$ and $\beta_C$ 
similar to the ones used in Figure \ref{eigenplot} (e.g. $\beta_I,\beta_C)=\{(-1.e-2,1.e-2),(-1.e-2,1.e-1)$), and the source
term is given by $f=\sin(x+iy)$. For problem \eqref{Max} (resp.
\eqref{lap}), the N\'ed\'elec (resp. $P_1$-Lagrange) finite elements are
used.  Figure~\ref{presults} exhibits convergence of the residual
for the  iterative procedures.
Clearly the case  $(\beta_I,\beta_C) =(0,0)$ is less performant. 

 More advanced numerical analysis of the proposed scheme will be the subject of
 forthcoming work. We shall in particular numerically and theoretically discuss
 optimal choices for the impedance parameters $\beta_C$ and $\beta_I$ and
 explore the possibility of using parameters that are operators of an
 appropriate negative
 order that would provide an asymptotic limit of  $|(\mathcal{T})^{m}_{n}|$ strictly
 smaller than $1$ (this is clearly not the case of constant parameters as
 indicated by Proposition \ref{asymptotics}).  
\section*{Acknowledgements}
Y. Boubendir gratefully acknowledges support from NSF through grant No. DMS-1319720.
\bibliographystyle{plain}
\bibliography{ref}

\begin{thebibliography}{10}

\bibitem{abramowitz1964handbook}
Milton Abramowitz and Irene~A Stegun.
\newblock {\em Handbook of mathematical functions: with formulas, graphs, and
  mathematical tables}.
\newblock Number~55. Courier Corporation, 1964.

\bibitem{MR3343738}
Ana Alonso~Rodr{\'{\i}}guez, Enrico Bertolazzi, Riccardo Ghiloni, and Alberto
  Valli.
\newblock Finite element simulation of eddy current problems using magnetic
  scalar potentials.
\newblock {\em J. Comput. Phys.}, 294:503--523, 2015.

\bibitem{MR2680968}
Ana Alonso~Rodr{\'{\i}}guez and Alberto Valli.
\newblock {\em Eddy current approximation of {M}axwell equations}, volume~4 of
  {\em MS\&A. Modeling, Simulation and Applications}.
\newblock Springer-Verlag Italia, Milan, 2010.
\newblock Theory, algorithms and applications.

\bibitem{MR3349682}
Ana Alonso~Rodr{\'{\i}}guez and Alberto Valli.
\newblock Finite element potentials.
\newblock {\em Appl. Numer. Math.}, 95:2--14, 2015.

\bibitem{boubendir2007analysis}
Yassine Boubendir.
\newblock An analysis of the bem-fem non-overlapping domain decomposition
  method for a scattering problem.
\newblock {\em Journal of computational and applied mathematics},
  204(2):282--291, 2007.

\bibitem{boubendir2012quasi}
Yassine Boubendir, Xavier Antoine, and Christophe Geuzaine.
\newblock A quasi-optimal non-overlapping domain decomposition algorithm for
  the helmholtz equation.
\newblock {\em Journal of Computational Physics}, 231(2):262--280, 2012.

\bibitem{gander2002optimized}
Martin~J Gander, Fr{\'e}d{\'e}ric Magoules, and Fr{\'e}d{\'e}ric Nataf.
\newblock Optimized schwarz methods without overlap for the helmholtz equation.
\newblock {\em SIAM Journal on Scientific Computing}, 24(1):38--60, 2002.

\bibitem{quarteroni1999domain}
Alfio Quarteroni and Alberto Valli.
\newblock {\em Domain decomposition methods for partial differential
  equations}.
\newblock Number CMCS-BOOK-2009-019. Oxford University Press, 1999.

\bibitem{NUMNUM20060}
Ana~Alonso Rodríguez, Ralf Hiptmair, and Alberto Valli.
\newblock A hybrid formulation of eddy current problems.
\newblock {\em Numerical Methods for Partial Differential Equations},
  21(4):742--763, 2005.

\bibitem{toselli2005domain}
Andrea Toselli and Olof Widlund.
\newblock {\em Domain decomposition methods: algorithms and theory}, volume~3.
\newblock Springer, 2005.

\end{thebibliography}
\end{document}